# Mathematics of Knowledge Refinement:

# Probabilistic Arithmetic, with no unknowns and no infinity.

## Part I. Generalized Probabilistic Arithmetic. Basic definitions and properties.

**Mikhail Luboschinsky**

An approach to build Probabilistic Arithmetic in which initial values of all correlated random variables are known, but with varying degrees of accuracy. As a result of the proposed Probabilistic Arithmetic operations, variable values, degrees of their accuracy and correlations are refined.

Probabilistic Generalized Addition (PGA) and Probabilistic Generalized Multiplication (PGM) operations on correlated random variables are defined and their basic properties identified and described:

- Proposed PGA and PGM operations possess inverse operations – subtraction and division.
- There is no difference between direct and inverse operations: addition and subtraction, multiplication and division (this is why these operations are called "Generalized").
- Division by approximately zero is possible and the result never equals to ∞, making this approach promising in computational and computer mathematics.
- Unlike the usual hyperbola, which, when the argument is changing from + ∞ to 0, has a gap at 0, the Generalized Probabilistic Hyperbola, under certain combination of initial accuracies, is continuous at 0. First, as the argument changes from +∞ to 0, it behaves like a typical hyperbola, monotonically increasing. However, after reaching certain maximum value, it starts to decrease to 0 monotonically.

Last property gives a hope that the proposed approach might be promising in Quantum Physics, as it would allow a more adequate description of macro and micro physics, as well as of the transition from one to the other.

==============================================================

**Математика уточнения знания:**

**Вероятностнаная математика без неизвестных и без бесконечности.**

**Часть I. Обобщённая Вероятностная арифметика. Основные определения и свойства.**

Михаил Любощинский

Предложен подход к построению Арифметики Случайных Чисел, при котором все исходные коррелированные случайные операнды считаются известными, но с разной степенью точности. В результате предложенных Вероятностных Арифметических Операций происходит уточнение значений, точностей и корреляций всех операндов.

Определены операции Обобщённого Вероятностного Сложения (ОВС) и Обобщённого Вероятностного Умножения (ОВУ) кореллированных случайных чисел, выявлены и описаны их основные свойства:

- Предложенные операции обладают обратными операциями.

- Нет разницы между прямыми и "обратными" операциями: между сложением и вычитанием, умножением и делением (почему эти операции и названы Обобщёнными).

- Возможно деление на приблизительно ноль и, в результате деления, никогда не возникает ∞, что делает этот подход очень перспективным в вычислительной и компьютерной математике.

- В отличие от обычной гиперболы, которая, при изменении аргумента от +∞ до -∞, претерпевает разрыв в точке 0, Обобщённая Вероятностная Гипербола, вначале, при изменении аргумента от +∞ к 0, ведёт себя как обычная гипербола, монотонно возрастая, но, затем, достигнув некоторого максимального значения, при некоторых соотношениях точностей, либо также имеет разрыв в окрестности 0, либо начинает монотонно убывать до 0, переходя в область отрицательных значений, сохраняя непрерывность.

Последнее свойство позволяет надеяться, что предлагаемый подход может оказаться перспективным в Квантовой Физике, т.к. позволит более адекватно описывать физику макро и микромира, и процессы перехода из макро в микро мир.

==========================================================

В работе [1] излагается метод решения плохо обусловленных систем линейных алгебраических уравнений

(1) $\boldsymbol{A\bar{x}} = \boldsymbol{\bar{f}},$

когда известны параметры распределения случайных матрицы $\boldsymbol{A}$ и вектора $\boldsymbol{\bar{f}}$. Кроме того, предполагается известной априорная информация об искомом векторе $\boldsymbol{\bar{x}}$. В цитируемой работе рассматривается случай, когда матрица $\boldsymbol{A}$ и векторы $\boldsymbol{\bar{f}},$ и $\boldsymbol{\bar{x}},$ принадлежат нормально распределенным совокупностям.

Суть метода состоит в том, что ищется максимум совместной плотности распределения, представляющей собою произведение плотностей распределения вероятностей матрицы $\boldsymbol{A}$, векторов $\boldsymbol{\bar{f}},$ и $\boldsymbol{\bar{x}},$ и невязки $\boldsymbol{h} = (\boldsymbol{\bar{f}} - \boldsymbol{A\bar{x}})^2$, в предположении, что задана дисперсия $\boldsymbol{\sigma^2}$ невязки $\boldsymbol{h}$.

Возникла идея попытаться построить, на основе этого подхода, Арифметику Случайных Чисел.

**<u>Постановка задачи.</u>**

Рассмотрим совокупность трех случайных чисел: $(\boldsymbol{\tilde{x}}, \boldsymbol{\tilde{y}}, \boldsymbol{\tilde{z}}),$ заданных совместной плотностью нормального распределения:

(2) $f(x, y, z, a, b, c, \Xi) = C \cdot \exp\left[\frac{-1}{2}(x-a,\ y-b,\ z-c)\Xi(x-a,\ y-b,\ z-c)^T\right]$

где $\boldsymbol{a}$, $\boldsymbol{b}$, $\boldsymbol{c}$- математические ожидания случайных чисел $\boldsymbol{\tilde{x}},$ $\boldsymbol{\tilde{y}}$ и $\boldsymbol{\tilde{z}}$, соответственно**,**

(3) $\Xi = \Sigma\Sigma^T = \begin{bmatrix} A & E & Z \\ E & B & H \\ Z & H & \Gamma \end{bmatrix}$ - матрица "точности" совместного распределения чисел $\boldsymbol{\tilde{x}}, \boldsymbol{\tilde{y}}, \boldsymbol{\tilde{z}},$ а $\Sigma = \Lambda^{-1}$, где $\Lambda$ - ковариационная матрица распределения.

Определим, содержательно, операции обобщённого вероятностного сложения или умножения, как уточнение функции совместной плотности распределения вероятности (2), за счёт использования дополнительной информации о том, что случайные операнды $\boldsymbol{\tilde{x}},$ $\boldsymbol{\tilde{y}}$ и $\boldsymbol{\tilde{z}}$ связаны соотношением $\boldsymbol{\tilde{x}} + \boldsymbol{\tilde{y}} = \boldsymbol{\tilde{z}}$ или $\boldsymbol{\tilde{x}} \cdot \boldsymbol{\tilde{y}} = \boldsymbol{\tilde{z}}$, т.е., предположим, что задана дисперсия $\theta^2$ невязки $\boldsymbol{\Delta_\oplus} = \boldsymbol{x} + \boldsymbol{y} - \boldsymbol{z},$ для сложения, или $\boldsymbol{\Delta_\odot} = \boldsymbol{x} \cdot \boldsymbol{y} - \boldsymbol{z},$ для умножения.

Для этого рассмотрим плотности распределения невязки $\Delta$:

(4) $f_{\Delta_\oplus}(x,y,z,\Theta) = C\cdot\exp\left[\frac{-\Theta}{2}(x+y-z)^2\right]$ - для сложения,

(5) $f_{\Delta_\odot}(x,y,z,\Theta) = C\cdot\exp\left[\frac{-\Theta}{2}(x\cdot y-z)^2\right]$ - для умножения.

Где $\Theta = \frac{1}{\theta^2}$ - точность операции сложения или умножения, а $\theta^2$ – дисперсия операции.

Введём следующие обозначения:

- $\tilde{\boldsymbol{x}} = {}_{\pm\alpha}\tilde{\boldsymbol{a}}$ , или, для краткости: $\tilde{\boldsymbol{x}} = {}_{\alpha}\tilde{\boldsymbol{a}}$ - случайное число $x$ со средним значением $a$ и дисперсией $\alpha^2$.
- $\left[{}_\alpha\tilde{a}, {}_\beta\tilde{b}\,{}_\theta\tilde{\oplus}\,{}_\gamma\tilde{c}\right]_\Xi = \left[\tilde{x},\tilde{y}\,{}_\theta\tilde{\oplus}\,\tilde{z}\right]_\Xi = [\tilde{x}\tilde{+}\tilde{y}\ {}_\theta\tilde{\cong}\ \tilde{z}]_\Xi = \left[{}_\alpha\tilde{a}\tilde{+}{}_\beta\tilde{b}\ {}_\theta\tilde{\cong}\ {}_\gamma\tilde{c}\right]_\Xi$ - обобщённое вероятностное сложение случайных чисел $\tilde{\boldsymbol{x}},\tilde{\boldsymbol{y}},\tilde{\boldsymbol{z}}$ с дисперсией операции сложения равной $\theta^2$, и матрицей точности $\boldsymbol{\Xi}$.
- $\left[{}_\alpha\tilde{a}, {}_\beta\tilde{b}\,{}_\theta\tilde{\odot}\,{}_\gamma\tilde{c}\right]_\Xi = \left[\tilde{x},\tilde{y}\,{}_\theta\tilde{\odot}\,\tilde{z}\right]_\Xi = [\tilde{x}\,\tilde{\cdot}\,\tilde{y}\ {}_\theta\tilde{\cong}\ \tilde{z}]_\Xi = \left[{}_\alpha\tilde{a}\,\tilde{\cdot}\,{}_\beta\tilde{b}\ {}_\theta\tilde{\cong}\ {}_\gamma\tilde{c}\right]_\Xi$ - обобщённое вероятностное умножение случайных чисел $\tilde{\boldsymbol{x}},\tilde{\boldsymbol{y}},\tilde{\boldsymbol{z}}$ с дисперсией операции умножения равной $\theta^2$, и матрицей точности $\boldsymbol{\Xi}$.

Тогда запись ${}_3\tilde{5}$ читается так: случайное число со средним значением равным 5 и средним квадратическим отклонением, равным 3, или, что то же: приблизительно 5, с дисперсией $3^2$.

***<u>Определение 1</u>***: Пусть тройка случайных чисел $(\tilde{\boldsymbol{x}},\tilde{\boldsymbol{y}},\tilde{\boldsymbol{z}})$ задана совместной плотностью распределения $f(x,y,z,a,b,c,\Xi)$, и допустимая точность операции сложения задана плотностью распределения невязки $f_{\Delta_\oplus}(x,y,z,\Theta)$.

Результатом операции Обобщенного Вероятностного Сложения (ОВС) тройки случайных чисел $(\tilde{\boldsymbol{x}},\tilde{\boldsymbol{y}},\tilde{\boldsymbol{z}})$ с дисперсией операции $\theta^2$ и матрицей точности $\boldsymbol{\Xi}$:

(6) $[\tilde{\boldsymbol{x}}\tilde{+}\tilde{\boldsymbol{y}}\ {}_\theta\tilde{\cong}\ \tilde{\boldsymbol{z}}]_\Xi = (\tilde{\boldsymbol{x}}',\tilde{\boldsymbol{y}}',\tilde{\boldsymbol{z}}')_{\Xi'}$,

является тройка случайных чисел $(\tilde{\boldsymbol{x}}',\tilde{\boldsymbol{y}}',\tilde{\boldsymbol{z}}')$, с совместной плотностью распределения $f_\oplus(x,y,z,a',b',c',\Xi')$, являющейся локальной аппроксимацией произведения плотностей распределения $f(x,y,z,a,b,c,\Xi)$ и $f_{\Delta_\oplus}(x,y,z,\Theta)$ в окрестности точки максимума:

(7) $f_\oplus(x,y,z,a',b',c',\Xi') = f(x,y,z,a,b,c,\Xi)\cdot f_{\Delta_\oplus}(x,y,z,\Theta) =$

$= C\cdot\exp\left[\frac{1}{2}[(x-a,y-b,z-c)\Xi(x-a,y-b,z-c)^T + \Theta(x+y-z)^2]\right]$,

где параметры распределения $a',b',c',\Xi'$ определяются следующим образом:

• $a', b', c'$ - уточнённые средние значения – являются решением задачи:

(8) $$F_{\oplus}(x,y,z,\Xi,\Theta)=\frac{1}{2}[(x-a,y-b,z-c)\Xi(x-a,y-b,z-c)^T+$$

$$+\Theta(x+y-z)^2]\ \rightarrow\ \min_{x,y,z}$$

• $\Xi'$ - уточнённые значения матрицы точности $\Xi$ для ОВС определяется из:

(9) $$\Xi'=\begin{bmatrix}\frac{\partial^2}{\partial x^2} & \frac{\partial^2}{\partial x\partial y} & \frac{\partial^2}{\partial x\partial z}\\ \frac{\partial^2}{\partial y\partial x} & \frac{\partial^2}{\partial y^2} & \frac{\partial^2}{\partial y\partial z}\\ \frac{\partial^2}{\partial z\partial x} & \frac{\partial^2}{\partial z\partial y} & \frac{\partial^2}{\partial z^2}\end{bmatrix}\cdot F_{\oplus}(x,y,z,\Xi,H)=\Xi+\Theta\begin{bmatrix}1 & 1 & -1\\ 1 & 1 & -1\\ -1 & -1 & 1\end{bmatrix}.$$

Из необходимых условий экстремума для (8), имеем:

(10) $$\nabla_{x,y,z}F_{\oplus}(x,y,z,\Xi,\Theta)=\left[\Xi+\Theta\begin{bmatrix}1 & 1 & -1\\ 1 & 1 & -1\\ -1 & -1 & 1\end{bmatrix}\right]\cdot\begin{bmatrix}x\\ y\\ z\end{bmatrix}-\Xi\cdot\begin{bmatrix}a\\ b\\ c\end{bmatrix}=0$$

Определяем координаты точки оптимума $\boldsymbol{x}^*$, $\boldsymbol{y}^*$ и $\boldsymbol{z}^*$ из (10):

(11) $$\begin{bmatrix}x^*\\ y^*\\ z^*\end{bmatrix}=\left[\Xi+\Theta\begin{bmatrix}1 & 1 & -1\\ 1 & 1 & -1\\ -1 & -1 & 1\end{bmatrix}\right]^{-1}\cdot\Xi\cdot\begin{bmatrix}a\\ b\\ c\end{bmatrix}$$

Таким образом, из (9) и (11) имеем, для точки оптимума $(a', b', c')$:

(12) $$\begin{bmatrix}a'\\ b'\\ c'\end{bmatrix}=\begin{bmatrix}x^*\\ y^*\\ z^*\end{bmatrix}=(\Xi')^{-1}\cdot\Xi\cdot\begin{bmatrix}a\\ b\\ c\end{bmatrix}$$

***Определение 2***: Пусть тройка случайных чисел $(\tilde{\boldsymbol{x}},\tilde{\boldsymbol{y}},\tilde{\boldsymbol{z}})$ задана совместной плотностью распределения $f(x,y,z,a,b,c,\Xi)$, и допустимая точность операции умножения задана плотностью распределения невязки $f_{\Delta_{\odot}}(x,y,z,\Theta)$.

Результатом операции Обобщенного Вероятностного Умножения (ОВУ) тройки случайных чисел $(\tilde{\boldsymbol{x}},\tilde{\boldsymbol{y}},\tilde{\boldsymbol{z}})$ с дисперсией операции $\theta^2$ и матрицей точности $\boldsymbol{\Xi}$:

(13) $$[\tilde{\boldsymbol{x}}\,\tilde{\cdot}\,\tilde{\boldsymbol{y}}\ {}_{\theta}\!\!\cong\ \tilde{\boldsymbol{z}}]_{\Xi}=(\tilde{\boldsymbol{x}}',\tilde{\boldsymbol{y}}',\tilde{\boldsymbol{z}}')_{\Xi'},$$

является тройка случайных чисел $(\tilde{\boldsymbol{x}}',\tilde{\boldsymbol{y}}',\tilde{\boldsymbol{z}}')$ с совместной плотностью распределения $f_{\odot}(x,y,z,a',b',c',\Xi')$, являющейся локальной аппроксимацией произведения плотностей распределения $f(x,y,z,a,b,c,\Xi)$ и $f_{\Delta_{\odot}}(x,y,z,\Theta)$ в окрестности точки максимума:

(14) $$f_{\odot}(x,y,z,a',b',c',\Xi') = f(x,y,z,a,b,c,\Xi)\cdot f_{\Delta_{\odot}}(x,y,z,\Theta) =$$

$$= C\cdot\exp\left[\frac{1}{2}[(x-a,y-b,z-c)\Xi(x-a,y-b,z-c)^T + \Theta(x\cdot y - z)^2]\right],$$

где параметры распределения $a', b', c', \Xi'$ определяются следующим образом:

- $a', b', c'$ - уточнённые средние значения – это значения $\boldsymbol{x}, \boldsymbol{y}, \boldsymbol{z}$, минимизирующие функционал:

(15) $$F_{\odot}(x,y,z,\Xi,\Theta) = \frac{1}{2}[(x-a,y-b,z-c)\Xi(x-a,y-b,z-c)^T +$$

$$+\Theta(x\cdot y - z)^2] \rightarrow \min_{x,y,z}$$

- $\Xi'$ - уточнённые значения матрицы точности $\Xi$ для ОВУ находим из:

(16) $$\Xi' = \begin{bmatrix} \frac{\partial^2}{\partial x^2} & \frac{\partial^2}{\partial x\partial y} & \frac{\partial^2}{\partial x\partial z} \\ \frac{\partial^2}{\partial y\partial x} & \frac{\partial^2}{\partial y^2} & \frac{\partial^2}{\partial y\partial z} \\ \frac{\partial^2}{\partial z\partial x} & \frac{\partial^2}{\partial z\partial y} & \frac{\partial^2}{\partial z^2} \end{bmatrix}\cdot F_{\odot}(x,y,z,\Xi,\Theta) = \Xi + \Theta\begin{bmatrix} x^2 & 2xy-z & -y \\ 2xy-z & y^2 & -x \\ -y & -x & 1 \end{bmatrix}$$

Для определения уточнённых средних значений $a', b', c'$, запишем необходимые условия экстремума функционала (15):

(17) $$\nabla_{x,y,z}F_{\odot}(x,y,z,\Xi,\Theta) = \Xi\begin{bmatrix} x-a \\ y-b \\ z-c \end{bmatrix} + \Theta\begin{bmatrix} y(xy-z) \\ x(xy-z) \\ -(xy-z) \end{bmatrix} =$$

$$= \begin{bmatrix} A & E & Z \\ E & B & H \\ Z & H & \Gamma \end{bmatrix}\cdot\begin{bmatrix} x-a \\ y-b \\ z-z \end{bmatrix} + \Theta\begin{bmatrix} y(xy-z) \\ x(xy-z) \\ -(xy-z) \end{bmatrix} = 0.$$

К сожалению, пока, не удалось найти аналитического решения системы нелинейных уравнений (17).

Для поиска экстремума функционала (15) использовались численные методы.

**Замечание 1.** При высокой точности операции $\Theta = \frac{1}{\theta^2}$, , задачи (8) и (15), геометрически, можно интерпретировать, как задачи нахождения точки касания эллипсоида равного уровня, порождённого функцией (2) $f(x,y,z,a,b,c,\Xi)$, и плоскостью $x+y=z$ или гиперболоидом $x*y=z$.

Таким образом, результатом Вероятностных Операций Обобщённого Сложения и Умножения, определёных и для статистически зависимых операндов, являются уточнённые ожидаемые значения всех трёх исходных операндов, а не только одного неизвестного, как в обычной Арифметике. Уточнение значений исходных операндов происходит в количестве, обратно пропорциональном их исходных точностей: чем выше априорная точность операнда, тем меньше изменяется (уточняется) его априорное значение – больше уточняется то, что мы хуже знаем, но, принципиально, уточняются все операнды. Кроме того, результатом обобщённых операций является уточнённая матрица точности совместного распределения операндов, описывающая точности компонентов результата операции и их корреляционную связь.

**Замечание 2.** Последнее свойство Обобщённых Вероятностных Операций может стать основой для построения эффективной теории контроля точности в компьютерной математике – и в вопросах хранения информации и, собственно, последовательного контроля точности: от учёта точности исходной информации, до точности конечных результатов всех вычислений.

**Замечание 3.** При очень высокой точности «известных» операндов, и очень низкой точности «неизвестного», Обобщённые Операции дают результаты, с высокой точностью приближающиеся к «обычной» арифметике.

Особенность обобщенных операций состоит в том, что, собственно, результатом операции будет умножение или деление (сложение или вычитание), или, даже, что-то промежуточное между умножением и делением (сложением или вычитанием), в зависимости от исходных значений дисперсионной матрицы $\Lambda$ и исходных математических ожиданий операндов.

В Обобщённой Вероятностной Арифметике, таким образом, пропадает разница между «прямыми» и «обратными» операциями, между сложением и вычитанием, или умножением и делением: всё определяется соотношением априорных точностей исходных операндов операций. Именно поэтому, введённые операции названы обобщёнными.

Возможны следующие основные случаи:

1. **Априорная точность операнда $\tilde{z}$ меньше, чем операндов $\tilde{x}$ и $\tilde{y}$.**

Обобщённая операция будет, в основном, уточнением ожидаемого значения операнда $\tilde{\boldsymbol{z}}$**,** и, таким образом, операция будет сложением, или умножением, соответственно.

### 1.1 Численный пример Обобщённого Вероятностного Сложения («сложение»):

$$\left[{}_{1}\widetilde{\mathbf{1}}\tilde{+}{}_{5}\widetilde{\mathbf{10}}\ {}_{\theta}{\cong}\ {}_{10}\widetilde{\mathbf{50}}\right]_{\Xi},$$

где $\theta = 0.1$, $\quad \Xi = \Sigma\Sigma^T = (\Lambda\Lambda^T)^{-1} = \begin{bmatrix} 1^2 & 0 & 0 \\ 0 & 5^2 & 0 \\ 0 & 0 & 10^2 \end{bmatrix}^{-1} = \begin{bmatrix} 1 & 0 & 0 \\ 0 & 0.04 & 0 \\ 0 & 0 & 0.01 \end{bmatrix}.$

Содержательно, эта операция имеет следующий смысл: “Найти наиболее вероятные значения случайных чисел $\boldsymbol{x}$ со средним значением $\boldsymbol{a = 1}$ и дисперсией $\mathbf{1}$, $\boldsymbol{y}$, со средним значением $\boldsymbol{b = 10}$ и дисперсией $\mathbf{5^2}$, и $\boldsymbol{z}$, со средним значением $\boldsymbol{c = 50}$ и дисперсией $\mathbf{10^2}$**,** обеспечивающих дисперсию $\boldsymbol{\theta^2}$ случайной величины невязки $\boldsymbol{h = x + y - z}$ равной $\mathbf{10^{-2}}$”. (При этом точность операции $\Theta = \frac{1}{\theta^2} = 10^2$).

Так как точность $\boldsymbol{C = 10^{-2}}$ операнда $\boldsymbol{z}$ меньше, чем операндов $\boldsymbol{x}$ и $\boldsymbol{y}$, то, в этом случае, операция ближе к сложению, т.е., больше происходит уточнение $\boldsymbol{z}$, чем $\boldsymbol{x}$ и $\boldsymbol{y}$.

В соответствии с (9):

$$\Xi' = \Xi + \Theta \begin{bmatrix} 1 & 1 & -1 \\ 1 & 1 & -1 \\ -1 & -1 & 1 \end{bmatrix} = \begin{bmatrix} 1 & 0 & 0 \\ 0 & 0.04 & 0 \\ 0 & 0 & 0.01 \end{bmatrix} + 100 \begin{bmatrix} 1 & 1 & -1 \\ 1 & 1 & -1 \\ -1 & -1 & 1 \end{bmatrix} =$$

$$= \begin{bmatrix} 101 & 100 & -100 \\ 100 & 100.04 & -100 \\ -100 & -100 & 100.01 \end{bmatrix}.$$

В соответствии с (12):

$$\begin{bmatrix} a' \\ b' \\ c' \end{bmatrix} = (\Xi')^{-1} \cdot \Xi \cdot \begin{bmatrix} a \\ b \\ c \end{bmatrix} = \begin{bmatrix} 0.992 & -7.936 \cdot 10^{-3} & -7.936 \cdot 10^{-3} \\ -0.198 & 0.802 & 0.198 \\ 0.794 & 0.794 & 0.206 \end{bmatrix} \cdot \begin{bmatrix} 1 \\ 10 \\ 50 \end{bmatrix} = \begin{bmatrix} 1.3095 \\ 17.7375 \\ 19.0501 \end{bmatrix}.$$

Таким образом, операция Обобщённого Вероятностного Сложения, в данном примере, переводит тройку операндов **(1, 10, 50)** в тройку ($\mathbf{1.3095}$, $\mathbf{17.7375}$, $\mathbf{19.0501}$)**:**

$$\left[{}_1\tilde{1} \tilde{+} {}_5\widetilde{10}\ {}_{0.1}\cong\ {}_{10}\widetilde{50}\right] = \left(1.\widetilde{3095},\ 17.\widetilde{7375},\ 19.\widetilde{0501}\right)_{\Xi'},$$

Величина невязки $\boldsymbol{h}$, в данном примере:

$$\boldsymbol{h = a' + b' - c'} = -3.095 \cdot 10^{-3}$$

На **<u>Рис. 1</u>** представлен трехмерный график, иллюстрирующий нахождение результата операции Обобщённого Вероятностного Сложения:

$$\left[{}_{\mathbf{1}}\tilde{\boldsymbol{a}} \tilde{+} {}_{\mathbf{5}}\widetilde{\mathbf{10}}\ {}_{\mathbf{0.1}}\cong\ {}_{\mathbf{10}}\widetilde{\mathbf{50}}\right], \text{ при } \quad -100 < \boldsymbol{a} < 100:$$

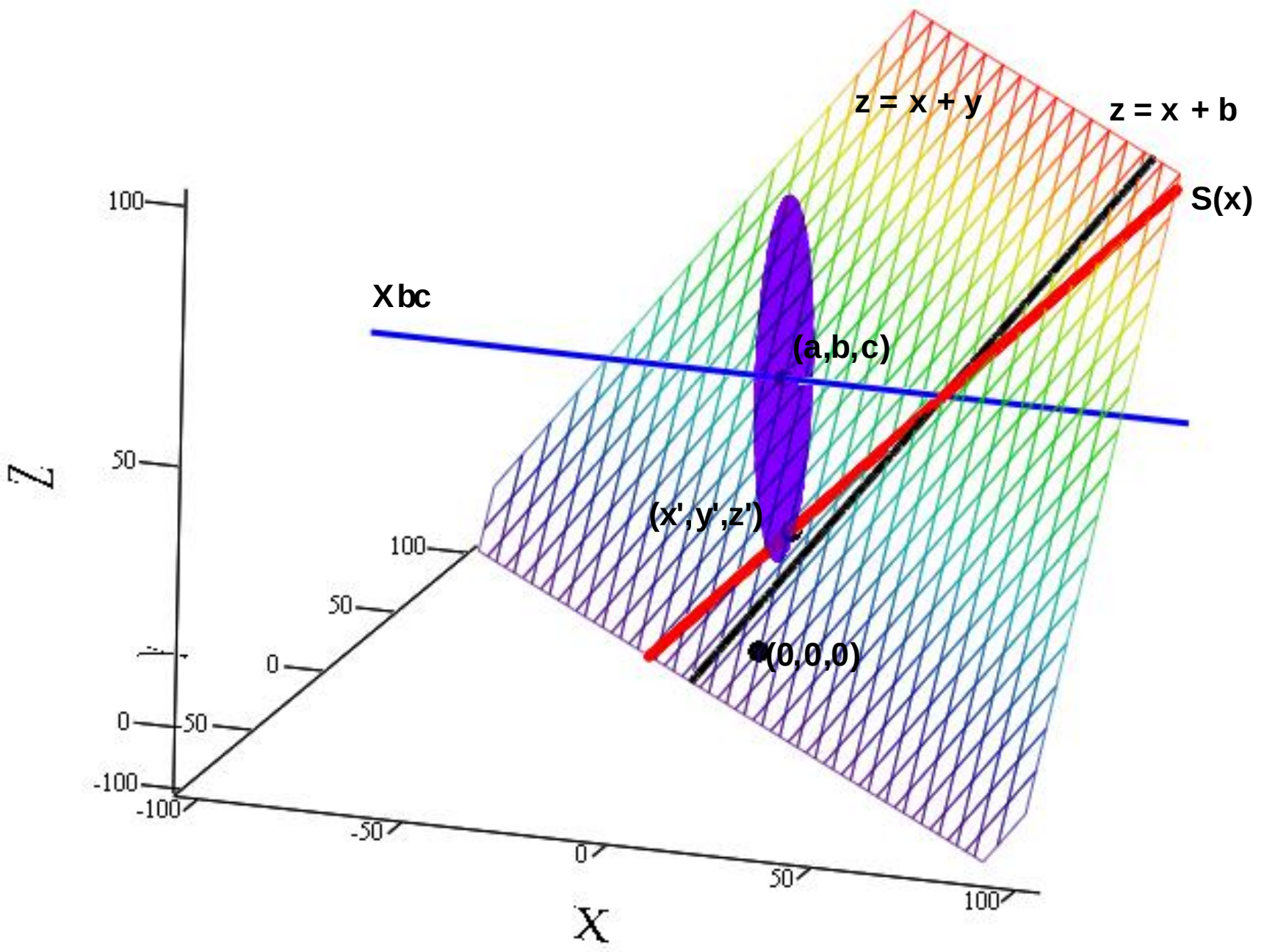

***Рис. 1. Построение результата операции Обобщённого Вероятностного Сложения.***

Обозначения:

- Точка **(0,0,0)** – начало координат.

- Точка **(a,b,c) = (1,10,50)** – центр распределения $\boldsymbol{f(x, y, z, a, b, c, \Xi)}$.

- Точка $(\boldsymbol{x', y', z'}) = (\mathbf{1.3095},\ \mathbf{17.7375},\ \mathbf{19.0501})$**:** – точка касания эллипсоидов равного уровня распределения функции $\boldsymbol{f(x, y, z, a, b, c, \Xi)}$ с плоскостью $\boldsymbol{z = x + y}$**.**

- **Xbc** – синяя прямая $\{\boldsymbol{x, y = b = 10, z = c = 50}\}$**,** по которой движется точка **(x,b,c)**, при изменении $\boldsymbol{x}$ от -100 до 100.

- $\boldsymbol{z = x + b}$ – чёрная прямая $\{\boldsymbol{x, y = 10, z = x + 10}\}$ (обычное сложение).

- $\boldsymbol{S(x)}$ – красная прямая, результат операции ОВС: геометрическое место точек $(\boldsymbol{x', y', z'})$ касания эллипсоида с центром в т. **(a,b,c),** порождённого функцией $\boldsymbol{f}$, при изменении координаты **a** от $-100$ до $+100$, т.е., при движении центра эллипсоида **(a,b,c)** по прямой **Xbc**.

### 1.2 Численный пример Обобщённого Вероятностного Умножения («умножение»):

$$\left[{}_{1}\widetilde{0.5}\ \tilde{\cdot}\ {}_{1}\widetilde{2}\ {}_{\theta}\cong\ {}_{10}\widetilde{5}\right]_{\Xi},$$

где $\theta = 0.1, \quad \Xi = \Sigma\Sigma^T = (\Lambda\Lambda^T)^{-1} = \begin{bmatrix} 1^2 & 0 & 0 \\ 0 & 1^2 & 0 \\ 0 & 0 & 10^2 \end{bmatrix}^{-1} = \begin{bmatrix} 1 & 0 & 0 \\ 0 & 1 & 0 \\ 0 & 0 & 0.01 \end{bmatrix}.$

Так как точность $\boldsymbol{C = 10^{-2}}$ операнда $\boldsymbol{z}$ меньше, чем операндов $\boldsymbol{x}$ и $\boldsymbol{y}$, то, в этом случае, операция ближе к умножению, т.е., происходит больше уточнение $\boldsymbol{z}$, чем $\boldsymbol{x}$ и $\boldsymbol{y}$.

Условия экстремума (17) имеют вид:

$$\begin{bmatrix} 1 & 0 & 0 \\ 0 & 1 & 0 \\ 0 & 0 & 0.01 \end{bmatrix} \cdot \begin{bmatrix} x - a \\ y - b \\ z - c \end{bmatrix} + 100 \cdot \begin{bmatrix} y(xy - z) \\ x(xy - z) \\ -(xy - z) \end{bmatrix} = 0$$

Решением этой системы нелинейных уравнений является тройка $(0.577,\ 2.022,\ 1.168)$.

Таким образом, $\left[{}_1\widetilde{0.5} \,\tilde{\cdot}\, {}_1\tilde{2} \;{}_{0.1}\!\cong\; {}_{10}\tilde{5}\right]_\Xi = (\widetilde{0.577}, \widetilde{2.022}, \widetilde{1.168})_{\Xi'}$ ,

Где, соответствии с (16):

$$\Xi' = \Xi + \Theta \begin{bmatrix} x^2 & 2xy - z & -y \\ 2xy - z & y^2 & -x \\ -y & -x & 1 \end{bmatrix} = 100 \cdot \begin{bmatrix} 0.343 & 1.167 & -2.022 \\ 1.167 & 4.099 & -0.577 \\ -2.022 & -0.577 & 1.000 \end{bmatrix}.$$

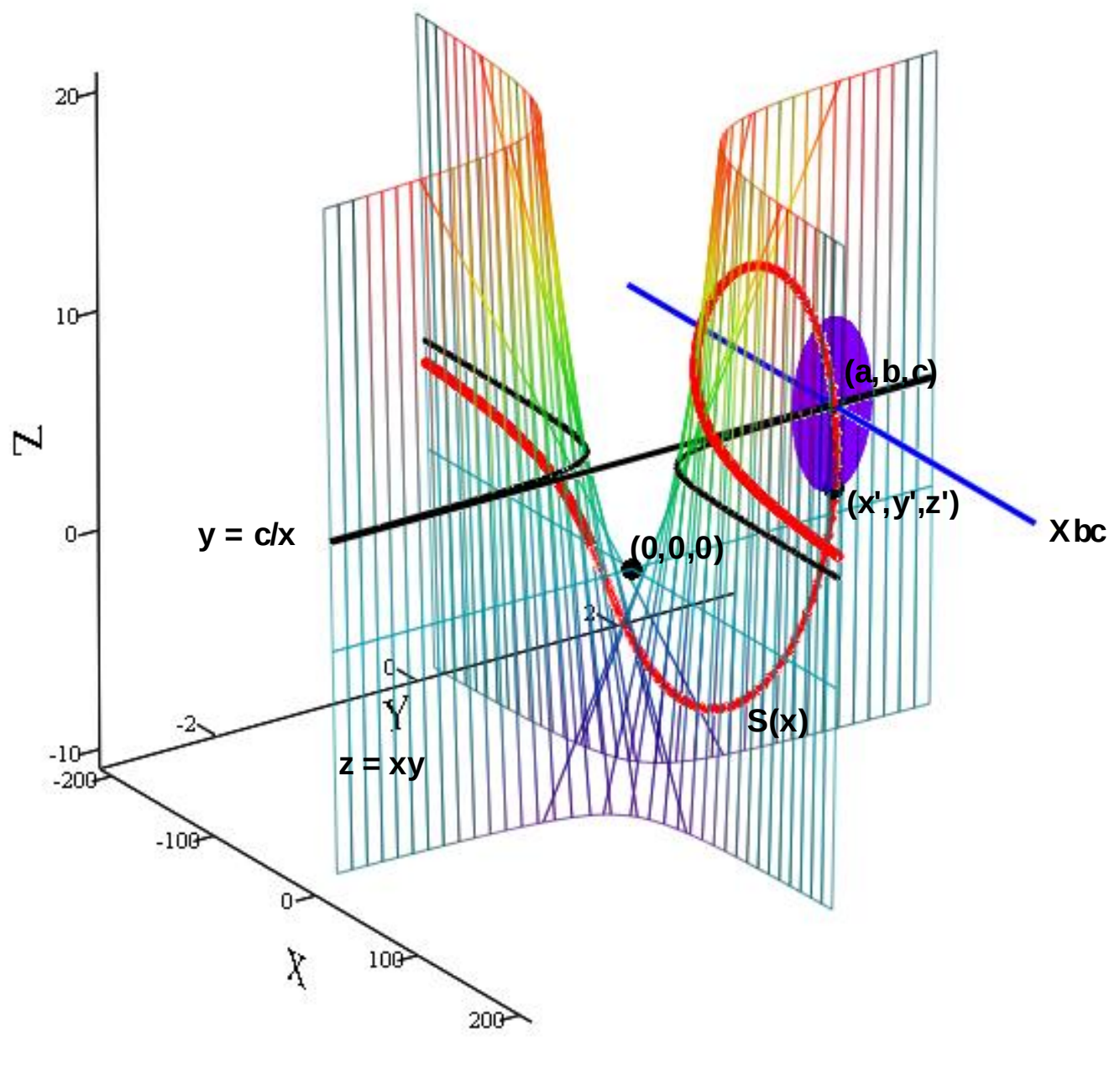

***<u>Рис. 2. Результат операции Обобщённого Вероятностного Умножения</u>***

$$\left[{}_{\mathbf{1}}\boldsymbol{\tilde{a}} \,\tilde{\cdot}\, {}_{\mathbf{1}}\boldsymbol{\tilde{2}} \;{}_{\mathbf{0.1}}\!\cong\; {}_{\mathbf{10}}\boldsymbol{\tilde{5}}\right]_{\boldsymbol{\Xi}}, \text{ при } -200 < \boldsymbol{a} < 200.$$

Обозначения:

• Точка $(\boldsymbol{x}', \boldsymbol{y}', \boldsymbol{z}') = (0.577,\ 2.022,\ 1.168)$ – точка касания эллипсоидов равного уровня распределения функции $\boldsymbol{f}(\boldsymbol{x}, \boldsymbol{y}, \boldsymbol{z}, \boldsymbol{a}, \boldsymbol{b}, \boldsymbol{c}, \boldsymbol{\Xi})$ с гиперболоидом $\boldsymbol{z} = \boldsymbol{x} \cdot \boldsymbol{y}$.

• $\mathbf{y} = \mathbf{c}/\mathbf{x}$ – чёрная кривая - обычная гипербола (обычное деление).

• $\boldsymbol{S}(\boldsymbol{x})$ – красная кривая - результат операции ОВУ: Обобщённая Вероятностная Гипербола: геометрическое место точек касания эллипсоида с центром в т. $(\boldsymbol{a}, \boldsymbol{b}, \boldsymbol{c})$, порождённого функцией $\boldsymbol{f}$, с гиперболоидом $\boldsymbol{z} = \boldsymbol{x} \cdot \boldsymbol{y}$, при изменении координаты $\boldsymbol{a} = \boldsymbol{x}$ от $-200\ до + 200$, т.е., при движении центра эллипсоида по прямой **Xbc**.

Хорошо видно первую особенность Вероятностной Гиперболы $\boldsymbol{S}(\boldsymbol{x})$ – это трехмерная спиралеобразная кривая, в отличие от обычной гиперболы $\boldsymbol{y} = \frac{\boldsymbol{c}}{\boldsymbol{x}}$, лежащей на плоскости $\boldsymbol{z} = \boldsymbol{c}$.

На **Рис. 3** представлен тот же трехмерный график, что и на **Рис. 2**, но в другом ракурсе - это вид со стороны оси $Y$ (проекция на плоскость $\{\boldsymbol{X}, \boldsymbol{Z}\}$):

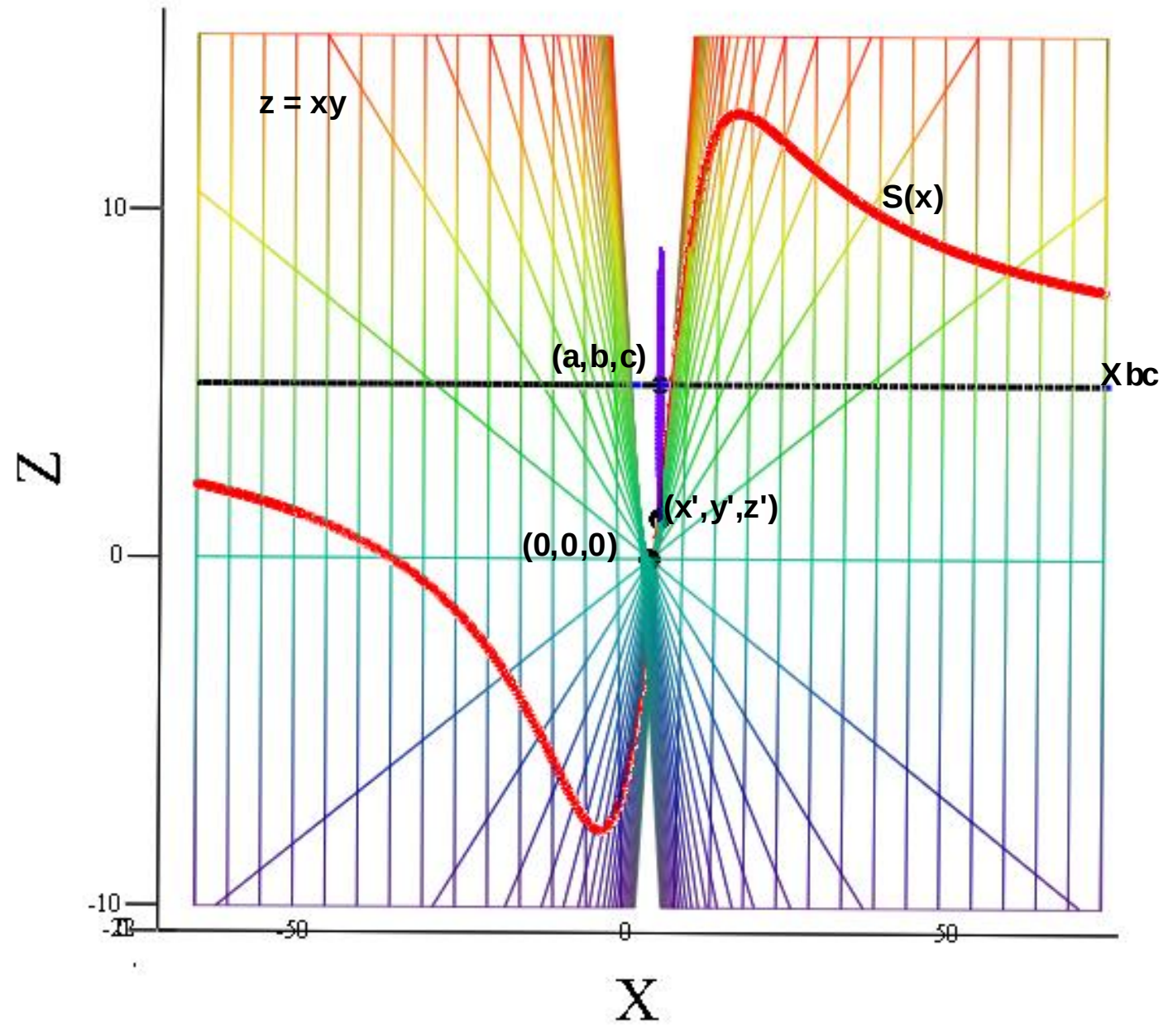

***Рис. 3. Обобщённое Вероятностное Умножение*** $\left[{}_{1}\tilde{\boldsymbol{a}}\ \tilde{\cdot}\ {}_{1}\tilde{\mathbf{2}}\ {}_{0.1}\widetilde{=}\ {}_{10}\tilde{\mathbf{5}}\right]_{\Xi}$

Проекция Обобщённой Вероятностной Гиперболы $\boldsymbol{S}(\boldsymbol{x})$ на плоскость $\{\boldsymbol{X}, \boldsymbol{Z}\}$,

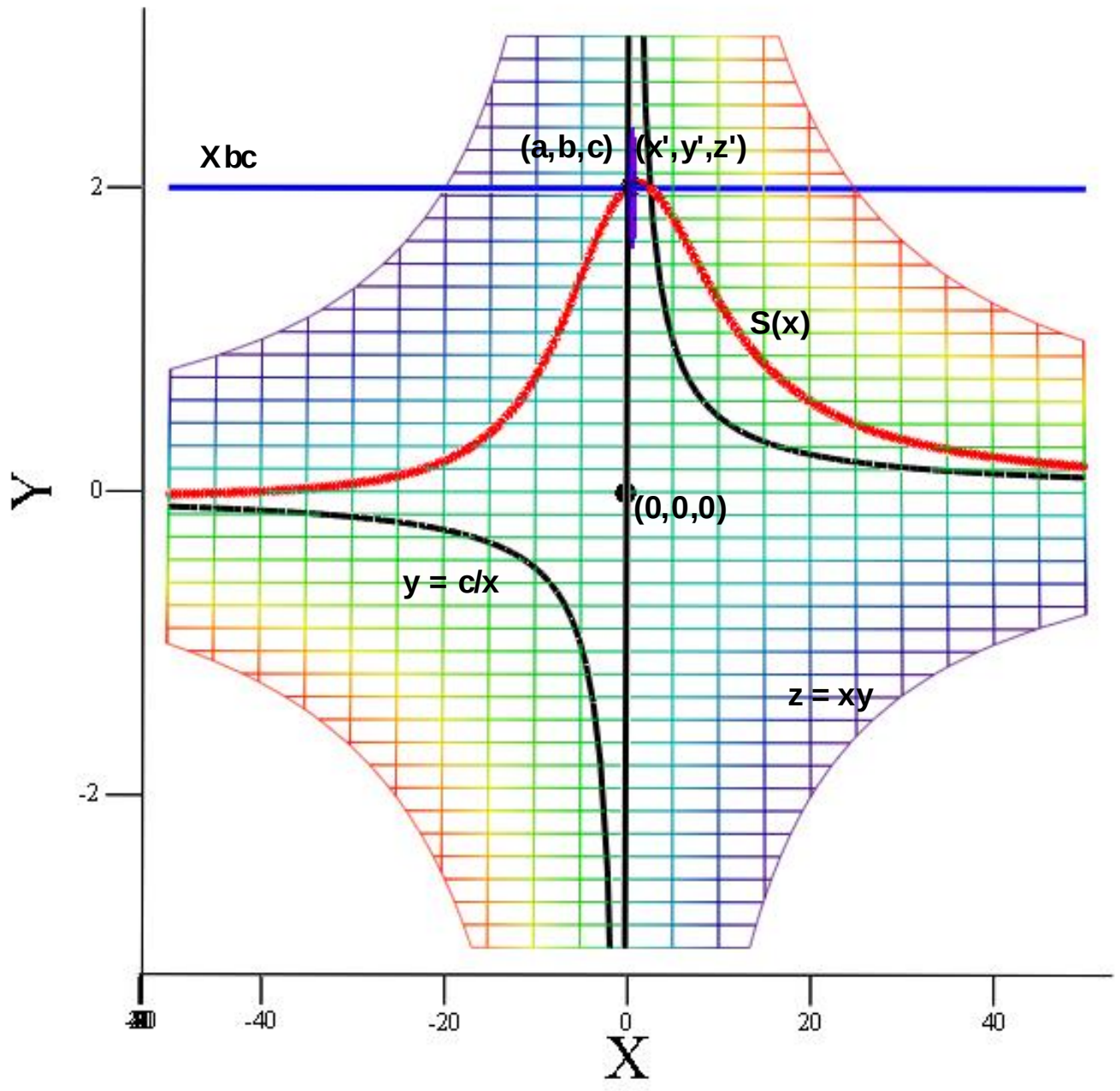

***Рис. 4. Обобщённое Вероятностное Умножение*** $\left[{}_{1}\tilde{a}\,\tilde{\cdot}\,{}_{1}\tilde{2}_{\,0.1}\cong\,{}_{10}\tilde{5}\right]_{\Xi}$

(проекция на плоскость $\{\boldsymbol{X},\boldsymbol{Y}\}$):

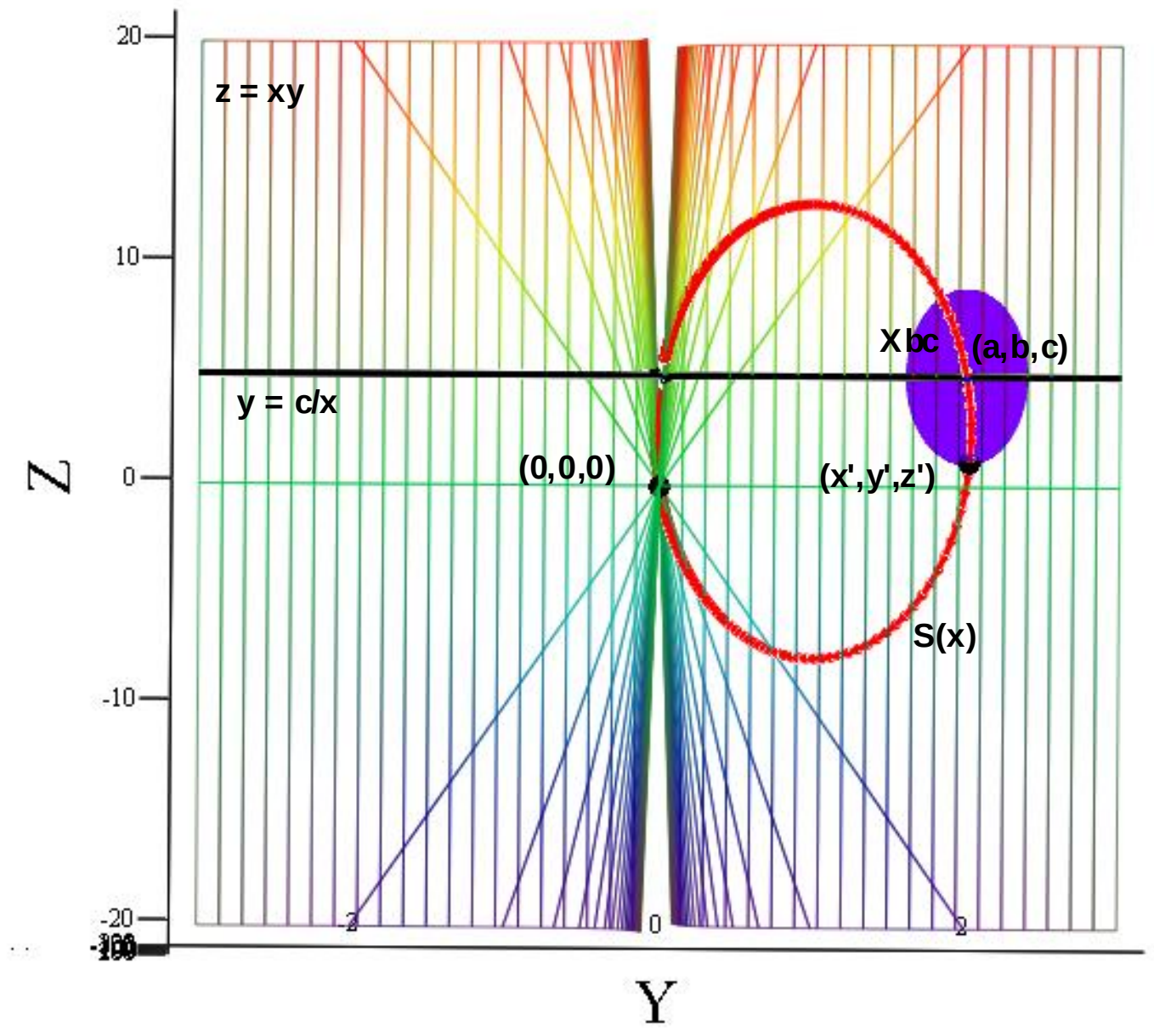

***Рис. 5. Обобщённое Вероятностное Умножение*** $\left[{}_{1}\tilde{a}\,\tilde{\cdot}\,{}_{1}\tilde{2}_{\,0.1}\cong\,{}_{10}\tilde{5}\right]_{\Xi}$

В этом ракурсе, проекция Обобщённой Вероятностной Гиперболы $\boldsymbol{S}(\boldsymbol{x})$ на плоскость $\{\boldsymbol{Y},\boldsymbol{Z}\}$, в общем случае, по-видимому, имеет форму эллипса. Этот вопрос требует специального рассмотрения. Возможно, некоторая замена переменных позволит найти аналитическое решение системы нелинейных уравнений (17).

**2. Априорная точность операнда $\tilde{x}$ (или $\tilde{y}$) меньше, чем операндов $\tilde{z}$ и $\tilde{y}$ (или $\tilde{z}$ и $\tilde{x}$, соответственно).**

Обобщённая операция, в этом случае, будет в основном, уточнением ожидаемого значения операнда $\tilde{\boldsymbol{x}}$ (или $\tilde{\boldsymbol{y}}$**,** соответственно)**,** и, таким образом, операция будет вычитанием или делением.

**2.1 Численный пример Обобщённого Вероятностного Сложения («вычитание»)**:

$$\left[{}_{1}\tilde{\boldsymbol{a}}\tilde{+}_{10}\tilde{\boldsymbol{2}}\;{}_{\theta}\tilde{\cong}\;{}_{3}\tilde{\boldsymbol{7}}\right]_{\Xi}\text{, где }\ \theta = 0.1,\quad \Xi = \begin{bmatrix} 1 & 0 & 0 \\ 0 & 0.01 & 0 \\ 0 & 0 & 0.111 \end{bmatrix}.$$

Так как дисперсия $10^2$ операнда $\boldsymbol{y}$ меньше, чем операндов $\boldsymbol{x}$ и $\boldsymbol{z}$, то, в этом случае, операция ближе к вычитанию, т.е., больше происходит уточнение $\boldsymbol{y}$, чем $\boldsymbol{x}$ и $\boldsymbol{z}$.

Операция Обобщённого Вероятностного Сложения («вычитания»), в данном примере, переводит тройку операндов **(1, 2, 7)** в тройку $(\tilde{\boldsymbol{x}}', \tilde{\boldsymbol{y}}', \tilde{\boldsymbol{z}}') = (\mathbf{1.036},\ \mathbf{5.636},\ \mathbf{6.672})$**:**

$$\left[{}_{1}\tilde{\boldsymbol{1}}\tilde{+}_{10}\tilde{\boldsymbol{2}}\;{}_{0.1}\tilde{\cong}\;{}_{3}\tilde{\boldsymbol{7}}\right] = (\mathbf{1.036},\ \mathbf{5.636},\ \mathbf{6.672})_{\Xi'},\quad \text{где }\ \Xi' = \begin{bmatrix} 101 & 100 & -100 \\ 100 & 100.01 & -100 \\ -100 & -100 & 100.111 \end{bmatrix}.$$

На **Рис. 6** представлен трехмерный график, иллюстрирующий нахождение результата операции Обобщённого Вероятностного Сложения («вычитания»):

$$\left[{}_{1}\tilde{\boldsymbol{a}}\tilde{+}_{10}\tilde{\boldsymbol{2}}\;{}_{10^2}\tilde{\cong}\;{}_{3}\tilde{\boldsymbol{7}}\right]\text{, при } -8 < \boldsymbol{a} < 8.$$

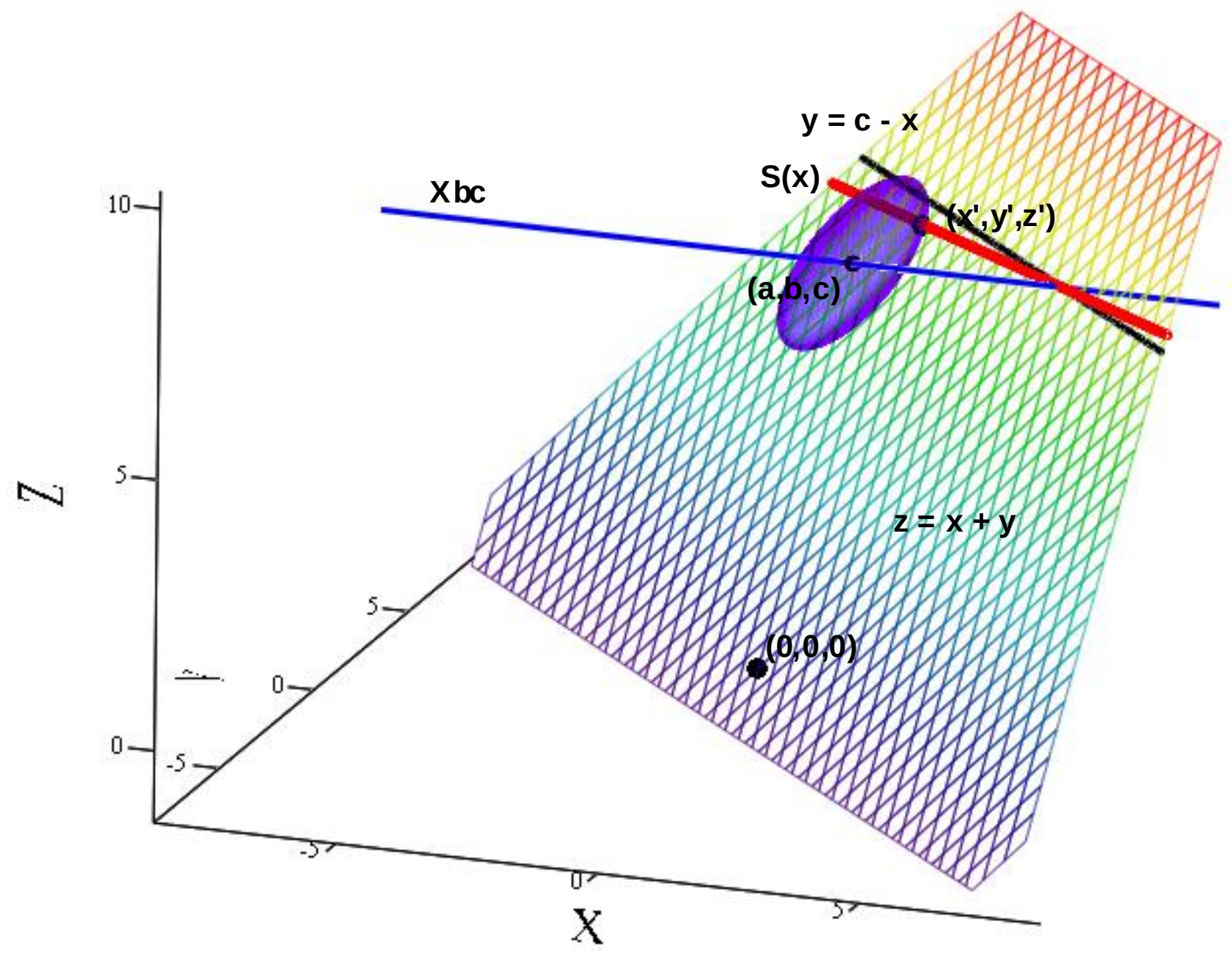

***Рис. 6. Обобщённое Вероятностное Сложение*** **(«вычитание»):** $\left[{}_{1}\widetilde{\boldsymbol{a}}\,\widetilde{\mp}\,{}_{10}\widetilde{\mathbf{2}}\;{}_{\mathbf{0.1}}\widetilde{\cong}\;{}_{\mathbf{3}}\widetilde{\mathbf{7}}\right]$

## 2.2 Численный пример Обобщённого Вероятностного Умножения («деления»):

$$\left[{}_{1}\widetilde{\mathbf{0.7}}\;\tilde{\cdot}\;{}_{10}\widetilde{\mathbf{2}}\;{}_{\theta}\widetilde{\cong}\;{}_{1}\widetilde{\mathbf{5}}\right]_{\boldsymbol{\Xi}}, \text{ где } \quad \theta = 0.1, \quad \Xi = \begin{bmatrix} 1 & 0 & 0 \\ 0 & 0.01 & 0 \\ 0 & 0 & 1 \end{bmatrix}.$$

Так как дисперсия $\mathbf{10^2}$ операнда $\boldsymbol{y}$ больше, чем операндов $\boldsymbol{x}$ и $\boldsymbol{z}$, то, в этом случае, операция ближе к делению, т.е., происходит больше уточнение $\boldsymbol{y}$, чем $\boldsymbol{x}$ и $\boldsymbol{z}$.

Решением системы нелинейных уравнений (17) является тройка $(0.908,\ 5.463,\ 4.962)$:

$$\left[{}_{1}\widetilde{\mathbf{0.7}}\;\tilde{\cdot}\;{}_{10}\widetilde{\mathbf{2}}\;{}_{\mathbf{0.1}}\widetilde{\cong}\;{}_{1}\widetilde{\mathbf{5}}\right]_{\boldsymbol{\Xi}} = (\widetilde{0.908}, \widetilde{5.463}, \widetilde{4.962})_{\Xi`}, \text{ где } \Xi' = 100\begin{bmatrix} 0.835 & 4.961 & -5.463 \\ 4.961 & 29.840 & -0.908 \\ -5.463 & -0.908 & 1.01 \end{bmatrix}.$$

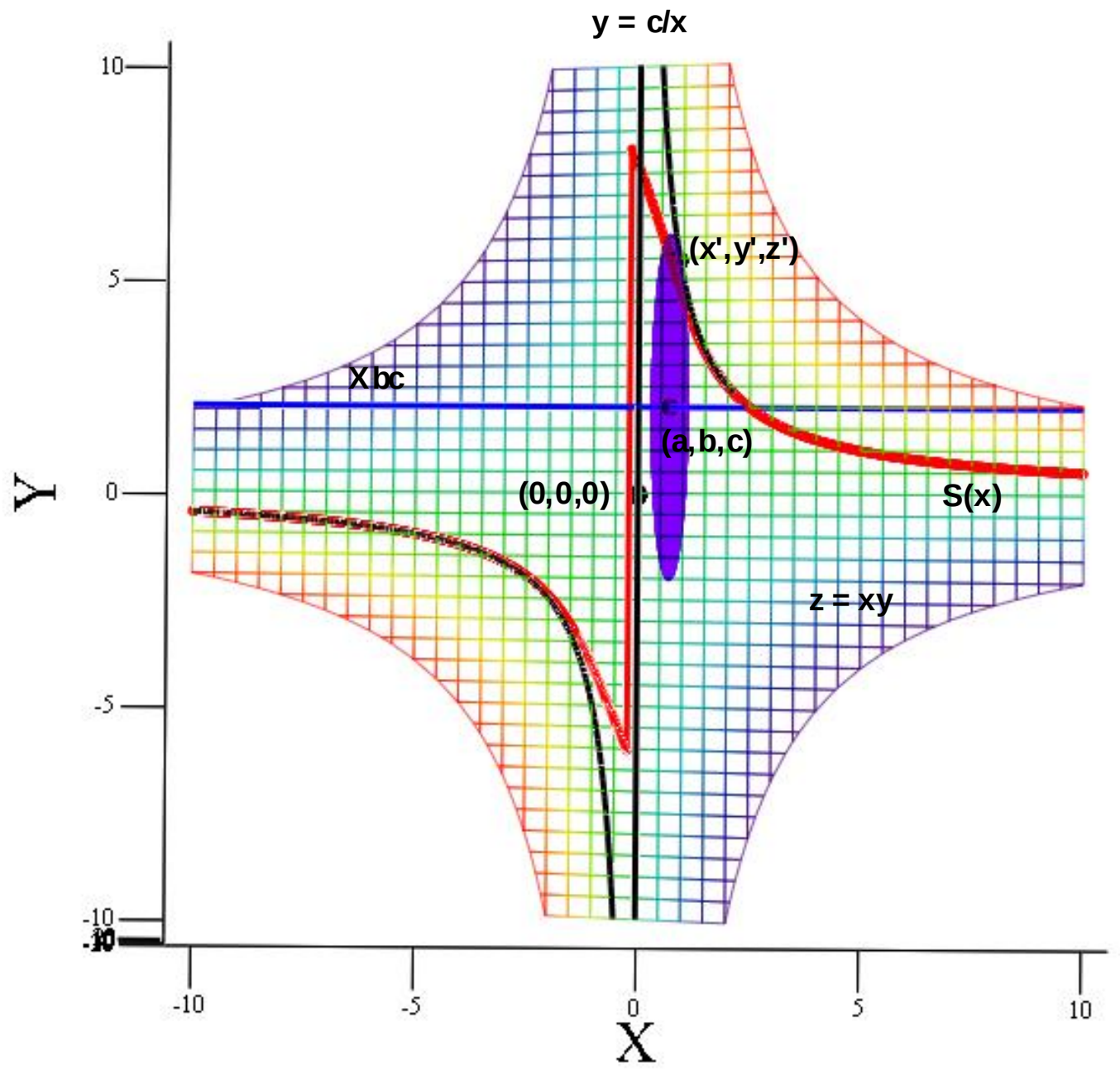

***Рис. 7. Обобщённое Вероятностное Умножение («деление»):*** $\left[{}_{1}\tilde{a}\,\tilde{\cdot}\,{}_{10}\tilde{2}\,{}_{0.1}\tilde{\cong}\;{}_{1}\tilde{5}\right]_{\Xi}$

Кривая $S(x)$ – Обобщённая Вероятностная Гипербола, при больших значениях $|x|$, хорошо совпадает с обычной гиперболой $y = \frac{c}{x}$, но затем, с уменьшением $|x|$, скорость роста, в сравнении с обычной гиперболой, уменьшается, и ордината $y$ достигает своего максимального значения. Характерная особенность Обобщённого Вероятностного умножения состоит в том, что, как видно из **Рис. 7**, при делении на приблизительно 0 не возникает ∞. В зависимости от исходных значений операндов и матрицы точности, при некотором значении $x$ (в окрестности 0), график Обобщённой Вероятностной Гиперболы $S(x)$ претерпевает разрыв, но соответствующие значения $y$ всегда конечны. (На **Рис. 7** заметен сдвиг в лево от 0 точки разрыва кривой $S(x)$, обусловленный тем, что ожидаемое значение операнда $y$ равно 2, а не 0. Это приводит к асимметричности Обобщённой Гиперболы относительно обеих осей $X$ и $Y$).

**Замечание 4.** Возможность деления на приблизительно 0 и конечность результата любой арифметической операции делает чрезвычайно перспективным применение Обобщённой Вероятностной Арифметики в вычислительной математике.

## 3. Априорная точность операнда $\tilde{z}$ больше точностей операндов $\tilde{x}$ и $\tilde{y}$.

В результате операции, происходит уточнение ожидаемых значений операндов и $\tilde{x}$, и $\tilde{y}$. В случае Обобщённого Сложения, операция представляет собой представление числа $\tilde{z}$ в виде суммы двух чисел $\tilde{\boldsymbol{x}}$ и $\tilde{\boldsymbol{y}}$. В случае Обобщённого Умножения, операция представляет собой представление числа $\tilde{\boldsymbol{z}}$ в виде произведения двух чисел $\tilde{\boldsymbol{x}}$ и $\tilde{\boldsymbol{y}}$ – факторизация числа $\tilde{\boldsymbol{z}}$.

### 3.1 Численный пример Обобщённого Вероятностного Умножения («факторизация»):

$$\left[{}_{\mathbf{10}}\tilde{\mathbf{1}} \mathbin{\dot{\tilde{}}} {}_{\mathbf{10}}\tilde{\mathbf{1}}\ {}_{\theta}\cong\ {}_{\mathbf{2}}\tilde{\mathbf{7}}\right]_{\Xi}, \text{ где } \theta = 0.01, \quad \Xi = \begin{bmatrix} 1 & 0 & 0 \\ 0 & 0.01 & 0 \\ 0 & 0 & 1 \end{bmatrix}.$$

Так как дисперсия $\mathbf{10^2}$ операндов $\boldsymbol{x}$ и $\boldsymbol{y}$ больше, чем операнда $\boldsymbol{z}$, то, в этом случае, происходит, больше уточнение одновременно и $\boldsymbol{x}$ и $\boldsymbol{y}$, чем операнда $\boldsymbol{z}$, т.е., факторизация числа $\boldsymbol{z}$.

Решением системы нелинейных уравнений (17) является тройка $(2.641,\ 2.641,\ 6.975)$:

$$\left[{}_{10}\tilde{1} \mathbin{\dot{\tilde{}}} {}_{10}\tilde{1}\ {}_{0.01}\cong\ {}_{2}\tilde{7}\right]_{\Xi} = (\widetilde{2.641}, \widetilde{2.641}, \widetilde{6.975})_{\Xi'} \text{ где } \Xi' = 10^4 \cdot \begin{bmatrix} 6.975 & 6.975 & -2.641 \\ 6.975 & 6.975 & -2.641 \\ -2.641 & -2.641 & 1 \end{bmatrix}.$$

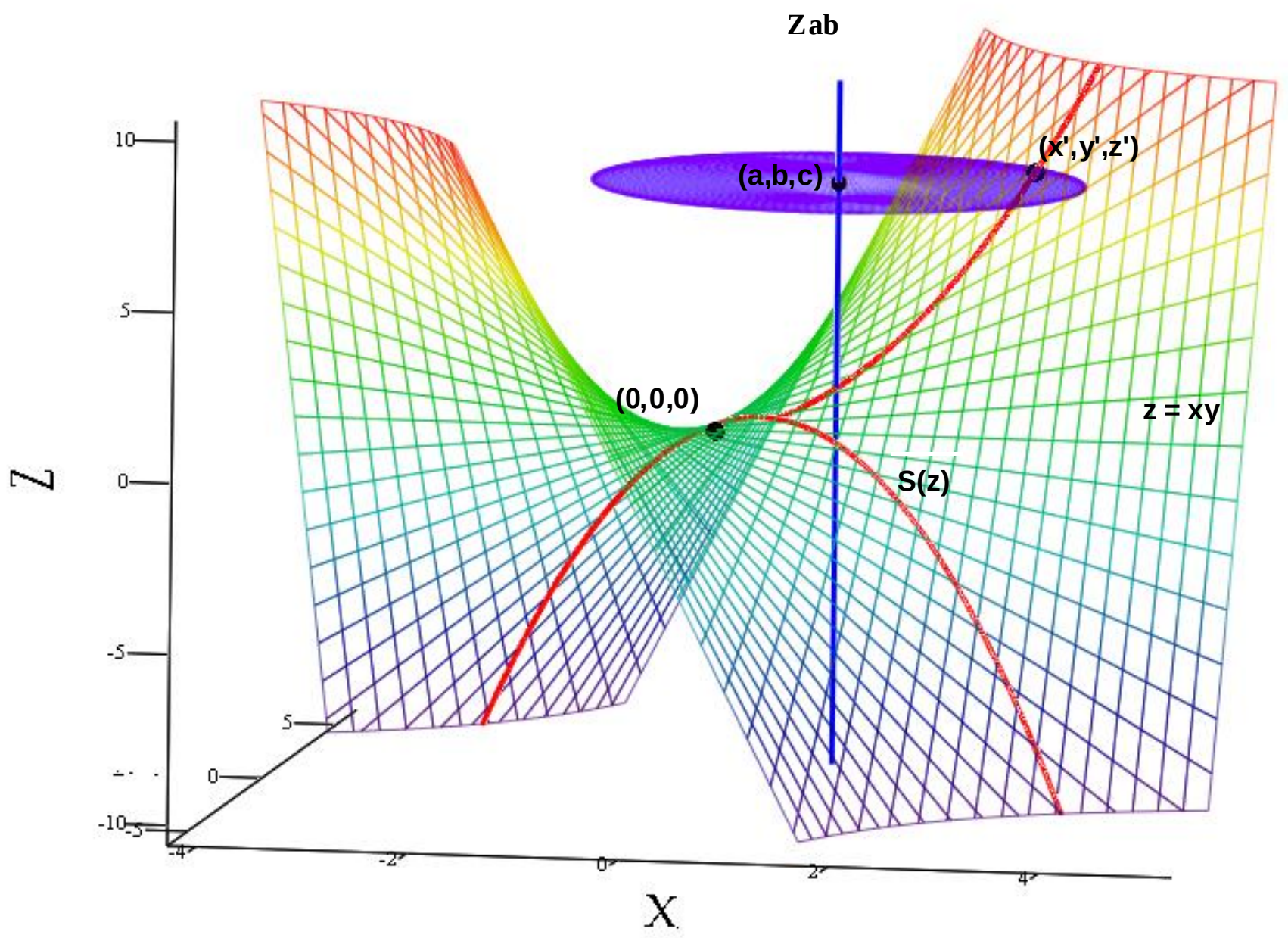

***Рис. 8. Обобщённое Вероятностное Умножение («факторизация»).***

$$\left[{}_{\mathbf{10}}\tilde{\mathbf{1}} \mathbin{\dot{\tilde{}}} {}_{\mathbf{10}}\tilde{\mathbf{1}}\ {}_{\mathbf{0.01}}\cong\ {}_{\mathbf{2}}\tilde{\boldsymbol{c}}\right], \text{ при } -10 < \boldsymbol{c} < 10.$$

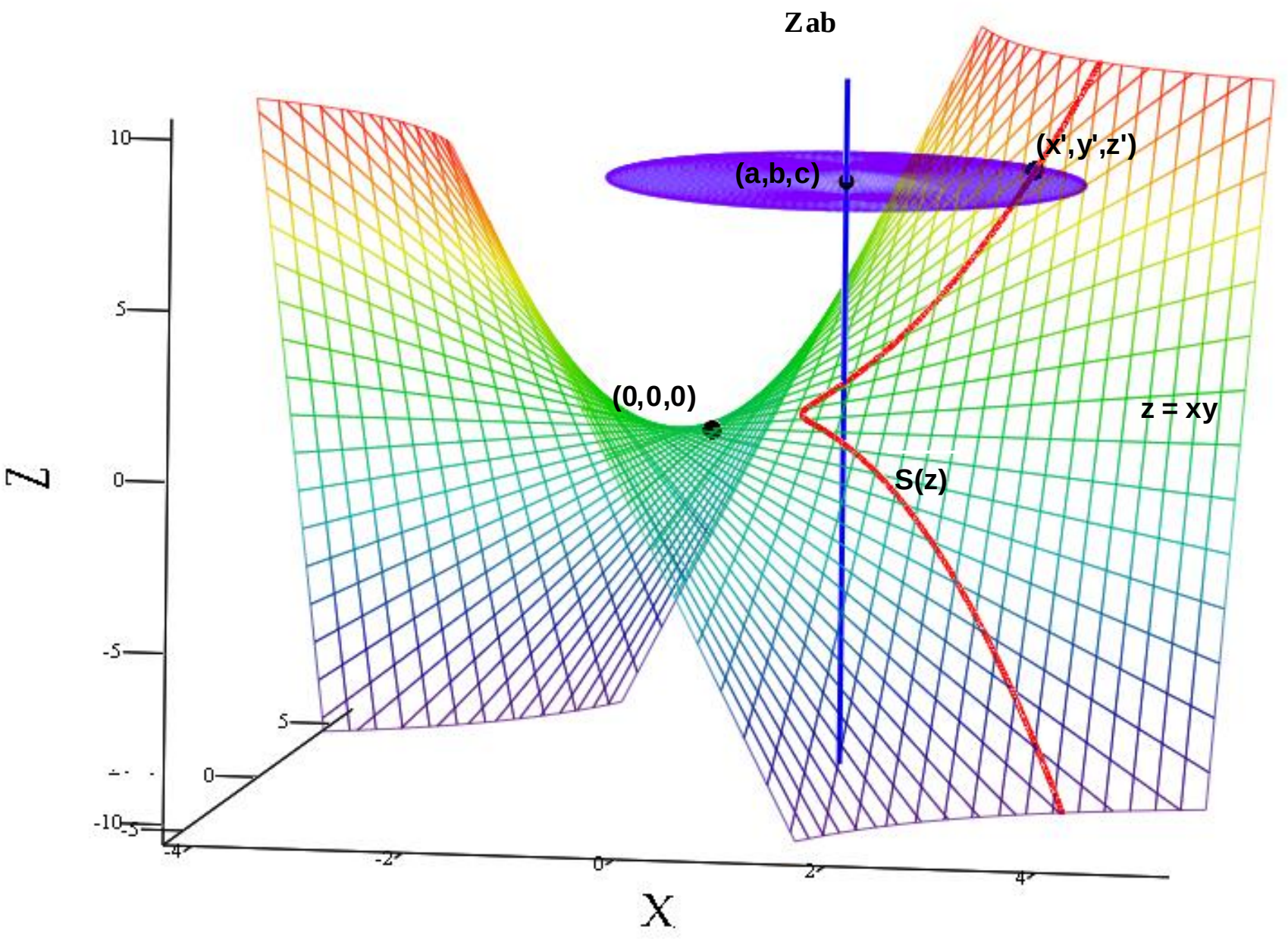

***<u>Рис. 9. Обобщённое Вероятностное Умножение («факторизация»).</u>***

$$\left[{}_{10}\widetilde{\mathbf{1.1}}\ \tilde{\cdot}\ {}_{10}\widetilde{\mathbf{1}}\ {}_{\mathbf{0.01}}\cong\ {}_{\mathbf{2}}\tilde{\boldsymbol{c}}\right],\ \text{при}\ -10 < \boldsymbol{c} < 10.$$

**4. <u>Априорная точность операнда $\tilde{x}$ (или $\tilde{y}$) больше точностей операндов $\tilde{z}$ и $\tilde{y}$ (или $\tilde{x}$, соответственно).</u>**

В результате операции, происходит, в основном, уточнение ожидаемых значений двух менее точных операндов $\tilde{\boldsymbol{z}}$ и $\tilde{\boldsymbol{y}}$ (или $\tilde{\boldsymbol{x}}$, соответственно), и операция будет чем-то промежуточным между сложением и вычитанием, или умножением и делением.

На **<u>Рис. 10</u>** представлен трехмерный график, иллюстрирующий нахождение результата операции Обобщённого Вероятностного Умножения:

$$\left[{}_{\mathbf{0.4}}\tilde{\boldsymbol{a}}\ \tilde{\cdot}\ {}_{\mathbf{10}}(-\tilde{\mathbf{2}})\ {}_{\mathbf{0.1}}\cong\ {}_{\mathbf{6}}\tilde{\mathbf{7}}\right],\text{при}\ -8 < \boldsymbol{a} < 8,$$

где $\theta = 0.1, \quad \varXi = \begin{bmatrix} 6.25 & 0 & 0 \\ 0 & 0.01 & 0 \\ 0 & 0 & 0.028 \end{bmatrix}.$

Так как точность операндов $\boldsymbol{y}$ и $\boldsymbol{z}$ меньше, чем операнда $\boldsymbol{x}$**,** то, в этом случае, происходит уточнение операндов $\boldsymbol{y}$ и $\boldsymbol{z}$, и операция представляет собой и умножение ( при значениях $\boldsymbol{x}$ близких к 0), и деление (при больших значениях $|\boldsymbol{x}|$**),** и нечто промежуточное между умножением и делением (при промежуточных значениях $|\boldsymbol{x}|$).

Операция Обобщённого Вероятностного Умножения, в данном примере, при $\boldsymbol{a} = \mathbf{1.2}$, переводит тройку операндов (1.2, -2, 7) в тройку $(\boldsymbol{x}', \boldsymbol{y}', \boldsymbol{z}') = \left(\widetilde{\mathbf{1.234}},\ \widetilde{\mathbf{4.205}},\ \widetilde{\mathbf{5.189}}\right)$:

$$\left[{}_{\mathbf{0.4}}\widetilde{\mathbf{2.3}}\ \tilde{\because}\ {}_{\mathbf{10}}(-\widetilde{\mathbf{2}})\ {}_{\mathbf{0.1}}\cong\ \widetilde{\mathbf{7}}\right]_{\Xi} = \left(\widetilde{\mathbf{1.234}},\ \widetilde{\mathbf{4.205}},\ \widetilde{\mathbf{5.189}}\right)_{\Xi'},$$

$$\text{где}\quad \Xi' = 10^2 \cdot \begin{bmatrix} 1.585 & 5.188 & -4.205 \\ 5.188 & 17.686 & -1.234 \\ -4.205 & -1.234 & 1 \end{bmatrix}$$

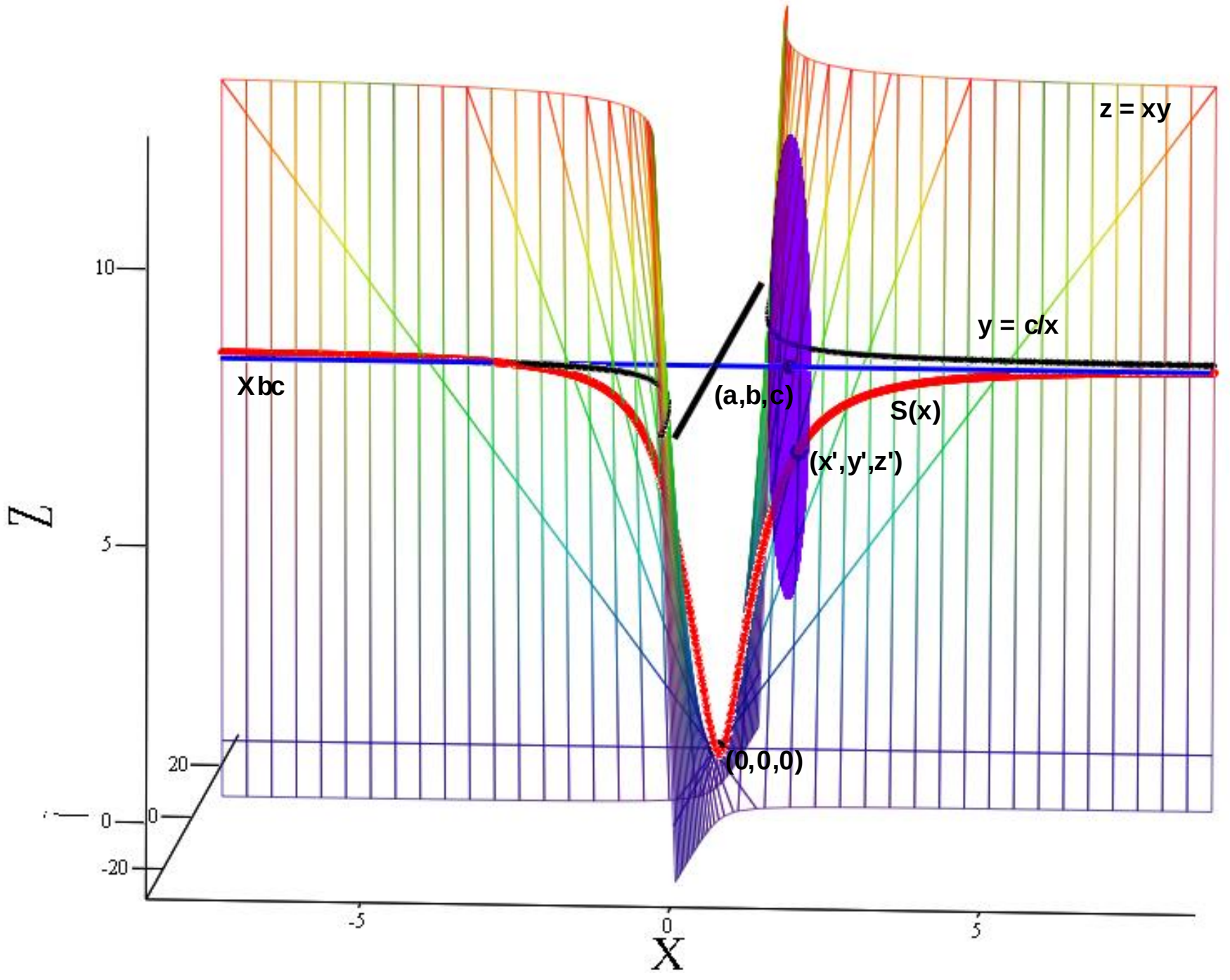

***Рис. 10. Обобщённое Вероятностное Умножение*** **(Умножение-Деление).**

$$\left[{}_{\mathbf{0.4}}\widetilde{\boldsymbol{a}}\ \tilde{\because}\ {}_{\mathbf{10}}(-\widetilde{\mathbf{2}})\ {}_{\mathbf{0.1}}\cong\ {}_{\mathbf{6}}\widetilde{\mathbf{7}}\right],\ \text{при}\ \ -8 < \boldsymbol{a} < 8.$$

Кривая **S(x)** – результат операции ОВУ: Обобщённая Вероятностная Гипербола: геометрическое место точек касания эллипсоида с центром в т. $(\boldsymbol{a}, \boldsymbol{b}, \boldsymbol{c})$, порождённого функцией $\boldsymbol{f}$, с гиперболоидом $\boldsymbol{z} = \boldsymbol{x} \cdot \boldsymbol{y}$, при изменении координаты $\boldsymbol{x} = \boldsymbol{a}$ от $-8\ до + 8$, т.е., при движении центра эллипсоида по прямой **Xbc**.

На **Рис. 11** представлен тот же трехмерный график, что и на **Рис. 10**, но в другом ракурсе - это вид сверху, со стороны оси $\boldsymbol{Z}$ (проекция на плоскость $\{\boldsymbol{X}, \boldsymbol{Y}\}$):

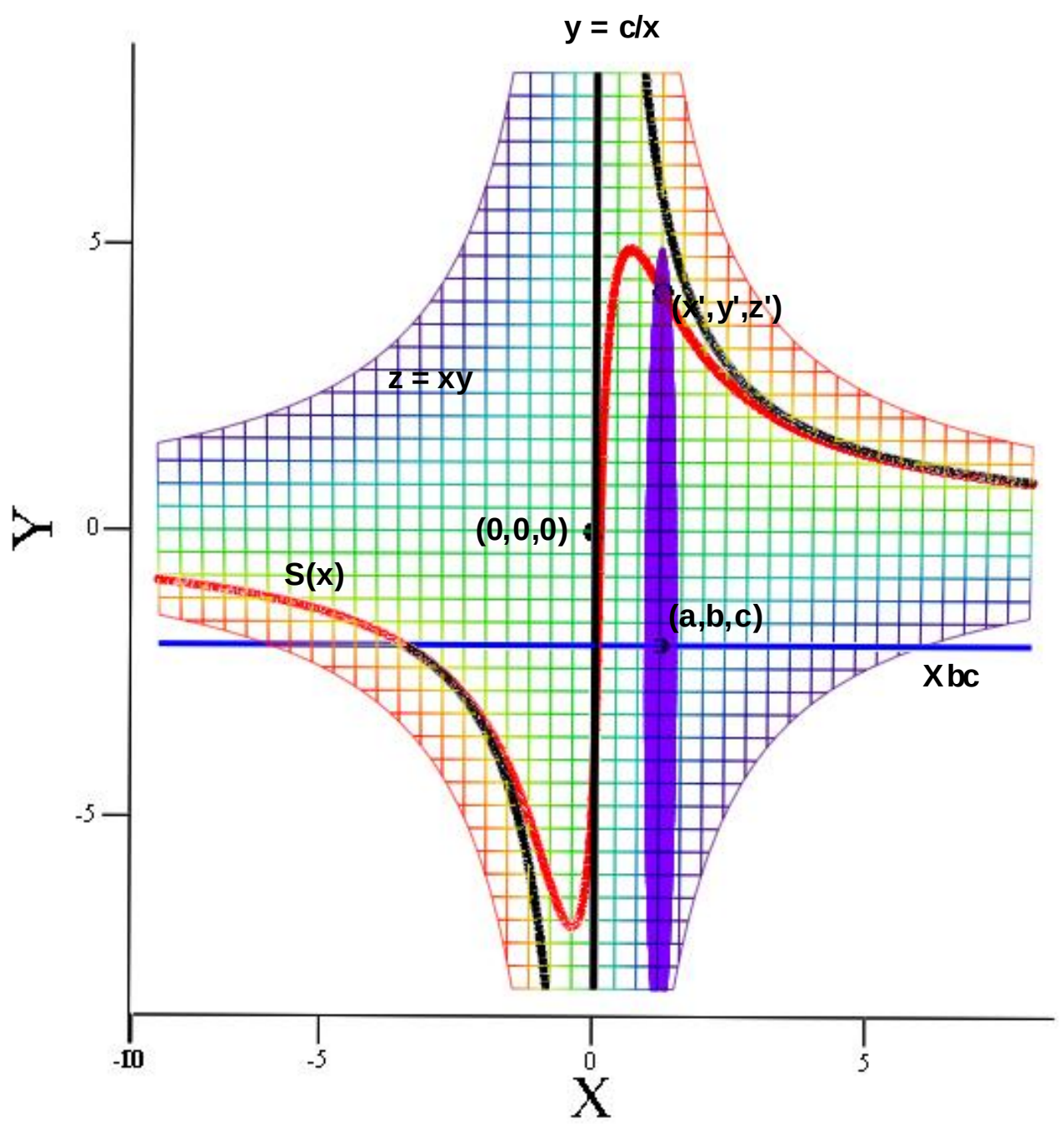

***Рис. 11. Обобщённое Вероятностное Умножение*** **(случай Умножения-Деления).**

На этом графике хорошо видна особенность Вероятностной Гиперболы $\boldsymbol{S(x)}$ в случае умножения-деления - в окрестности нуля, это умножение: ордината **y** возрастает с ростом **x**, достигает некоторое максимальное значение, и, затем, начинает уменьшаться, всё более приближаясь к обычной гиперболе, постепенно переходя в деление. Таким образом, в Обобщённой Вероятностной Арифметике возможно «деление» на приблизительно 0, и результат «деления» всегда конечен.

На **Рис. 12** - тот же график, что и на **Рис. 10 и 11**, но вид со стороны оси $\boldsymbol{X}$.

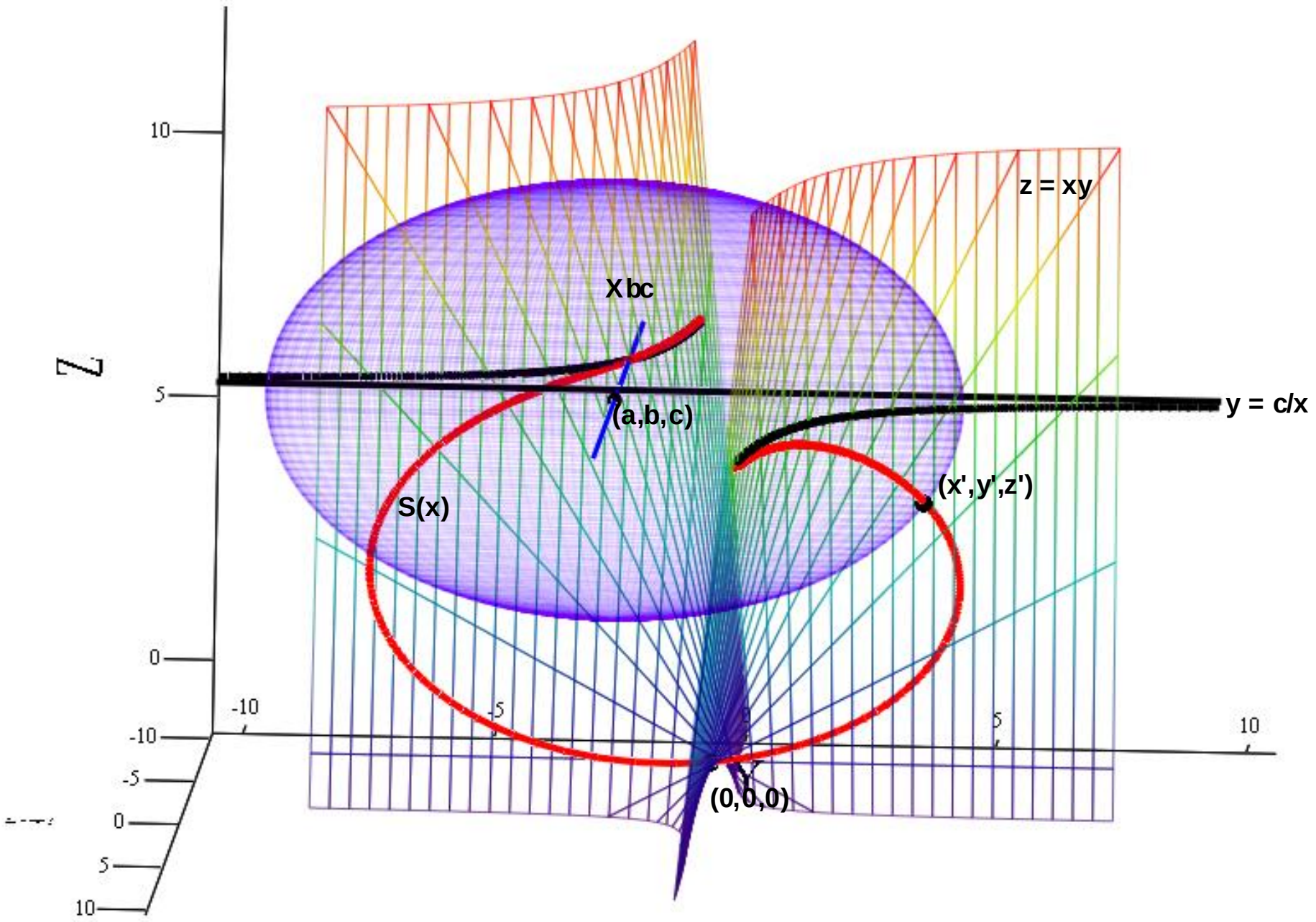

***Рис. 12. Обобщённое Вероятностное Умножение*** **(случай Умножения-Деления).**

Хорошо видно, что, в отличие от обычной гиперболы $\boldsymbol{y = \frac{c}{x}}$, лежащей на плоскости $\boldsymbol{z = c}$, Вероятностная Гипербола $\boldsymbol{S(x)}$ – это трехмерная спиралеобразная кривая, при больших значениях $|x|$ приближающаяся к обычной гиперболе, а в окрестности 0 существенно от неё отличающаяся. При $|x| \to 0$, Вероятностная Гипербола, вместо ухода в $\infty$, отходит от плоскости $\boldsymbol{z = c}$, и делает, в общем случае, петлю, проходя через 0. В частном случае (см. Рис.7) Вероятностная Гипербола претерпевает разрыв в 0, или в окрестности 0, в зависимости от исходных средних значений операндов, соотношений точностей и корреляций.

**Гипотеза.** То, что Обобщённая Вероятностная Гипербола при изменении $\boldsymbol{x}$ от больших значений к $\boldsymbol{0}$ совпадает с обычной гиперболой, а в окрестности $\boldsymbol{0}$ начинает вести себя очень своеобразно, позволяет высказать предположение, что предлагаемый подход позволит более адекватно описывать квантовую физику и микро и макро миров, и процессы перехода из одного мира в другой.

**<u>Библиография</u>**

1. Жуковский Е.Л. Метод наименьших квадратов для вырожденных и плохо обусловленных систем линейных алгебраических уравнений. Журнал Вычислительной Математики и Математической Физики. Том 17. №4, 1977.